\documentclass{article}

\usepackage[latin1]{inputenc}
\usepackage[english]{babel}
\usepackage{amsmath,amssymb}
\usepackage[dvips]{graphicx}
\usepackage{enumerate}
\usepackage{pstricks-add}
\usepackage{fullpage}
\usepackage{authblk} 

\newtheorem{Lem}{Lemma}
\newenvironment{lem}[1][]{\begin{Lem}\begin{itshape}\emph{#1 }}
{\end{itshape}\finenu\end{Lem}}
\newtheorem{Teo}[Lem]{Theorem}
\newenvironment{teo}[1][]{\begin{Teo}\begin{itshape}\emph{#1 }}
{\end{itshape}\finenu\end{Teo}}
\newtheorem{Cor}[Lem]{Corollary}
\newenvironment{cor}[1][]{\begin{Cor}\begin{itshape}\emph{#1 }}
{\end{itshape}\finenu\end{Cor}}
\newtheorem{Pro}[Lem]{Proposition}
\newenvironment{pro}[1][]{\begin{Pro}\begin{itshape}\emph{#1 }}
{\end{itshape}\finenu\end{Pro}}
\newtheorem{Defi}[Lem]{Definition}
\newenvironment{defi}[1][]{\begin{Defi}\begin{normalfont}\emph{#1 }}
{\end{normalfont}\fine\end{Defi}}
\newtheorem{Oss}[Lem]{Remark}
\newenvironment{oss}[1][]{\begin{Oss}\begin{normalfont}\emph{#1 }}
{\end{normalfont}\fineoss\end{Oss}}
\newtheorem{Es}[Lem]{Example}

\newenvironment{myproof}[1][]{\par\noindent\textbf{Proof{#1}. }}{\finedim\par}

\newcommand{\fine}{}
\newcommand{\finenu}{}
\newcommand{\finedim}{\hfill$\blacksquare$}
\newcommand{\fineoss}{\hfill$\square$}

\newcommand{\Banach}{\mathbb{R}^d}
\newcommand{\Closed}[2][]{\mathbb{K}}
\newcommand{\Fuzzy}[2][]{\mathbb{F}}
\newcommand{\salgebra}{\mathfrak{F}}
\newcommand{\borel}[1][]{\mathcal{B}_{#1}}
\newcommand{\prob}[1][]{\mathbb{P}_{#1}}
\newcommand{\misura}{\mu}
\newcommand{\dd}{{\ \rm d}}

\newcommand{\huku}[1]{{\Theta_{#1}}}
\newcommand{\detHuku}[1]{{\theta_{#1}}}
\newcommand{\hukuMin}[1]{{{H}_{#1}^{\bot}}}
\newcommand{\steiner}{\mathbf{ste} }
\newcommand{\genSteiner}{\mathbf{Ste} }
\newcommand{\unitSphere}{S^{\, d-1}}
\newcommand{\support}[1]{s_{#1}}
\newcommand{\supportBob}[1]{s^*_{#1}}

\newcommand{\Mink}{\oplus}
\newcommand{\Minkowski}{+}
\newcommand{\HukuDiff}{\ominus_{H}}
\newcommand{\HausDist}{\delta_H}
\newcommand{\indicator}[1]{\mathbb{I}_{#1}}
\newcommand{\racs}{{RaCS}}
\newcommand{\frv}{{FRV}}

\begin{document}

\title{A decomposition theorem for fuzzy set--valued random variables and a characterization of fuzzy random translation}
\author[1,2]{Giacomo Aletti\thanks{giacomo.aletti@unimi.it}}
\author[1]{Enea G. Bongiorno\thanks{enea.bongiorno@unimi.it}}
\affil[1]{Dipartimento di Matematica, Università degli studi di Milano}
\affil[2]{ADAMSS CENTRE (ADvanced Applied Mathematical and Statistical Sciences)}

\maketitle

\abstract{Let $X$ be a fuzzy set--valued random variable (\frv{}), and $\huku{X}$ the family of all fuzzy sets $B$ for which the Hukuhara difference $X\HukuDiff B$ exists $\mathbb{P}$--almost surely. In this paper, we prove that $X$ can be decomposed as $X(\omega)=C\Mink Y(\omega)$ for $\mathbb{P}$--almost every $\omega\in\Omega$, $C$ is the unique deterministic fuzzy set that minimizes $\mathbb{E}[d_2(X,B)^2]$ as $B$ is varying in $\huku{X}$, and $Y$ is a centered \frv{} (i.e. its generalized Steiner point is the origin). This decomposition allows us to characterize all \frv{} translation (i.e. $X(\omega) = M \Mink \indicator{\xi(\omega)}$ for some deterministic fuzzy convex set $M$ and some random element in $\Banach$). In particular, $X$ is a \frv{} translation if and only if the Aumann expectation $\mathbb{E}X$ is equal to $C$ up to a translation.
\\%
Examples, such as the Gaussian case, are provided.

\bigskip


\section*{Introduction}

It is widely known (e.g. \cite[Theorem~6.1.7]{li:ogu:kre02}) that a Gaussian fuzzy random variable may be decomposed as
\begin{equation}\label{eq:gauss_decomposition}
X=\mathbb{E} X \Mink \indicator{\xi },
\end{equation}
where $\mathbb{E} X$ is the expectation of $X$ in the Aumann sense, $\xi$ is a Gaussian random element in $\Banach$ with $\mathbb{E} \xi=0$ and $\indicator{A}:\Banach\to\{0,1\}$ denotes the indicator function of any $A\subseteq \Banach$
\[
\indicator{A}(x)=\left\{
\begin{array}{ll}
1, & \textrm{if } x\in A,\\
0, & \textrm{otherwise}.
\end{array}
\right.
\]
We write $\indicator{ a }$ instead of $\indicator{\{a\} }$ whenever $A=\{a\}$ is a singleton.
Roughly speaking, a Gaussian \frv{} $X$ is just a deterministic fuzzy set (its expected value $\mathbb{E} X$) up to a Gaussian translation $\xi$ which carries out all the randomness of $X$.
In this view, Equation~\eqref{eq:gauss_decomposition} entails a \emph{randomness defuzzification} for the Gaussian \frv{} $X$ according to which the underlying probability structure can be defined just only on $\Banach$ and no longer on $\Fuzzy[kc]{\Banach}$, the space of normal fuzzy sets with compact convex level sets.
Such randomness defuzzification occurs whenever a \frv{} $X$ is a random translation of a deterministic fuzzy set $M$. In this paper we provide a characterization for random translations by means of a suitable decomposition theorem that holds for any \frv{}.
In particular, given a centered \frv{} $X$, we define the family $\huku{X}$ of all deterministic $B\in \Fuzzy[kc]{\Banach}$ for which the Hukuhara difference $X\HukuDiff B$ exists almost surely. We show that this set is not empty, convex and closed in $(\Fuzzy[kc]{\Banach}, d_2)$, where $d_2$ corresponds to the $L^2$ metric in the space of support functions. Further,
\[
C = \underset{U\in\huku{X}}{\operatorname{arg\, min}}\, \mathbb{E}(d_2(X,U)^2)
\]
is unique and there exists a \frv{} $Y$ such that $X(\omega) =C\Mink Y(\omega)$; in some sense, $C$ and $Y$ are the deterministic part (with respect to $\Mink$) and the random part of $X$ respectively.
\\%
Since, the Aumann expectation $\mathbb{E}X$ is the (unique) Frèchet expectation with respect to $d_2$, i.e.
\[
\mathbb{E}X= \underset{U\in\Fuzzy[kc]{\Banach}}{\operatorname{arg\, min}}\, \mathbb{E}(d_2(X,U)^2),
\]
we obtain immediately that a \frv{} $X$ is a random translation of $C$ (i.e. $Y(\omega)$ is almost surely a singleton) if and only if $\mathbb{E}X$ is equal to $C$.

The paper is organized as follow. Section~\ref{sec: preliminary results} introduces necessary notations and literature results.  Section~\ref{sec:hukuhara_set} studies properties of the Hukuhara set $\huku{X}$ whilst Section~\ref{sec:hukuhara_decomposition} presents the decomposition theorem of \frv{} and the characterization of \frv{} translation.

\section{Preliminaries}
\label{sec: preliminary results}

Denote by $\Closed[kc]{\Banach}$ the class of non--empty compact convex subsets of $\Banach$, endowed with the Hausdorff metric
\[
\HausDist(A,B) = \max\{\sup_{a\in A} \inf_{b\in B}\|a-b\| ,
\sup_{b\in B} \inf_{a\in A}\|a-b\|\},%
\]
and the operations
\[
A\Minkowski B = \{a+b:a\in A,\ b\in B\},\qquad  \lambda\cdot A=\lambda A =\{\lambda a: a\in A\}\ \textrm{ with }\lambda> 0.
\]
For a non--empty closed convex set $A\subset \Banach$ the \emph{support function} $\support{A}:\unitSphere \to \mathbb{R}$ is defined by
\[
\support{A}(x)=\sup\{ \langle x,a \rangle  : a\in A \},\qquad \textrm{for } x\in\unitSphere,
\]
where $\langle \cdot,\cdot \rangle$ is the scalar product in $\Banach$ and $\unitSphere=\{x\in\Banach : \|x\|=1\}$ is the unit sphere in $\Banach$.
The \emph{Steiner point} of $A\in \Closed[kc]{\Banach}$ is defined by
\[
\steiner(A) = \frac{1}{v_d} \int_{\unitSphere} \support{A}(x)\, x \dd \lambda(x)
\]
where $x\in\unitSphere$ varies over the unit vectors of $\Banach$, $\lambda$ is the Lebesgue measure on $\unitSphere$, and $v_d$ is the volume of the unit ball of $\Banach$.

\paragraph{Fuzzy Sets.}
A {\em fuzzy set} is a map $\nu: \Banach \to [0,1]$. Let $\Fuzzy[kc]{\Banach}$ denote the family of all fuzzy sets $\nu$, which satisfy the following conditions.
\begin{enumerate}
\item
$\nu$ is an upper semicontinuous function, i.e.\ for each $\alpha\in (0,1]$, the \emph{cut set} or the \emph{$\alpha$--level set}    $\nu_\alpha=\{x\in\Banach : \nu (x)\ge \alpha \}$ is a closed subset of $\Banach$.

\item
$\nu$ is normal; i.e. $\nu_1=\{x\in\Banach : \nu (x)=1\}\neq\emptyset$.

\item
The support set $\nu_{0}=\overline{\{x\in\Banach : \nu(x)>0\}}$ of $\nu$ is compact; hence every $\nu_\alpha$ is compact for $\alpha\in(0,1]$.

\item
For any $\alpha\in [0,1]$, $\nu_\alpha$ is a convex subset of $\Banach$.
\end{enumerate}
For any $\nu\in\Fuzzy[kc]{\Banach}$ define the \emph{support function} of $\nu$ as follows:
\[
\support{\nu} (x,\alpha) =\left\{
\begin{array}{ll}
\support{\nu_\alpha}(x) & {\rm if } \ \alpha >0, \\
\support{\nu_{0}}(x) & {\rm if } \ \alpha =0,
\end{array}
\right.
\]
for $(x,\alpha)\in\unitSphere\times [0,1]$.
Let us endow $\Fuzzy[kc]{\Banach}$ with the operations
\[
(\nu^1\Mink\nu^2)_\alpha = \nu^1_\alpha \Minkowski \nu^2_\alpha,\qquad (\lambda\odot\nu^1)_\alpha = \lambda \cdot\nu^1_\alpha,\ \textrm{ with }\lambda> 0
\]
(so that $(\Fuzzy[kc]{\Banach},\Mink,\cdot)$ is a convex cone), and with the metrics
\begin{align*}
\HausDist^{\infty}(\nu^1, \nu^2)  & = \sup \{\alpha\in [0,1] :  \HausDist (\nu^1_\alpha, \nu^2_\alpha)\},
\\
d_2(\nu^1, \nu^2)  & = \left(\int_0^1 \int_{\unitSphere} | \support{\nu^1}(\alpha, u) - \support{\nu^2}(\alpha, u) |^2 \dd u \dd\alpha \right)^{\frac{1}{2}}.
\end{align*}
It is known that $(\Fuzzy[kc]{\Banach}, \HausDist^\infty)$ is a complete metric space while $(\Fuzzy[kc]{\Banach}, d_2)$ is not (cf. \cite[Chapter~7]{dia:klo94}).
The \emph{generalized Steiner point} of $A\in \Fuzzy[kc]{\Banach}$ is defined by
\[
\genSteiner(A) = \int_{[0,1]} \steiner(A_\alpha) \dd\alpha,
\]
where $\dd\alpha$ is the Lebesgue measure on $[0,1]$. In other words, $\genSteiner(A)$ may be seen as a weighted average of steiner points of the level sets of $A$. The following properties are satisfied (cf. \cite{vet:nav06}).
\begin{enumerate}
\item
For any $A\in\Fuzzy[kc]{\Banach}$, $\genSteiner(A)\in A_{0}$.
\item
For any $A,B\in\Fuzzy[kc]{\Banach}$, $\genSteiner(A\Mink B) = \genSteiner(A)+\genSteiner(B)$.

\item
$\genSteiner  : \Fuzzy[kc]{\Banach} \to \Banach$ is continuous.
\end{enumerate}

\paragraph{On the support function for fuzzy sets.}
It is known that the support function for a fuzzy set $\nu\in\Fuzzy[kc]{\Banach}$ can be defined equivalently on the closed unit ball $B(0,1)=\{x\in\Banach : \|x\|\le 1 \}\subset \Banach$ instead of the unit sphere $\unitSphere$ by
\[
\begin{array}{rccl}
\supportBob{\nu}:  & B(0,1) & \to & \mathbb{R} \\
 & x & \mapsto & \supportBob{\nu}(x)=\max \{ \langle x , y\rangle : y\in\Banach, \nu(y)\ge \|x\| \}.
\end{array}
\]
In particular, the following relationship between support function definitions hold
\begin{align*}
\forall (x,\alpha)\in \unitSphere \times [0,1], \qquad \support{\nu} (x,\alpha) & =
\left\{
\begin{array}{ll}
\supportBob{\nu}(\alpha x),  & \textrm{if } \alpha\neq 0;
\\
\sup_{y\in \nu_{0}} \langle y, x \rangle, & \textrm{if } \alpha= 0.
\end{array}
\right.
\\
\forall x\in B(0,1), \qquad \supportBob{\nu} (x) & =
\left\{
\begin{array}{ll}
\|x\| \, \support{\nu} \left( \frac{x}{\|x\|},\|x\| \right),  & \textrm{if } x\neq 0;
\\
0, & \textrm{if } x= 0.
\end{array}
\right.
\end{align*}
In \cite{bob85}, the author prove that a function $f:B(0,1)\to \mathbb{R}$ is a support function of some fuzzy set $\nu\in\Fuzzy[kc]{\Banach}$ if and only if the following six properties are satisfied:
\begin{enumerate}[(Property.1)]
\item
$f$ is upper semicontinuous, i.e.,
\[
f(x) = \limsup_{y\to x} f(y), \qquad \forall x\in B(0,1).
\]

\item
$f$ is positively semihomogeneous, i.e.,
\[
\lambda f(x) \le f(\lambda x), \qquad \forall\lambda\in (0,1], \forall x\in B(0,1).
\]

\item
$f$ is quasiadditive, i.e.,
\[
\|x\| f\left(\lambda \frac{x}{\|x\|}\right) \le \|x_1\| f\left(\lambda \frac{x_1}{\|x_1\|}\right) +
\|x_2\| f\left(\lambda \frac{x_2}{\|x_2\|}\right),
\]
for every $\lambda\in (0,1]$, and $x,x_1,x_2\in \Banach\setminus\{0\}$, with $x=x_1+x_2$.

\item
$f$ is normal, i.e.,
\[
f(x)+f(-x)\ge 0, \qquad \forall x\in B(0,1).
\]

\item
$f(\cdot)/\|\cdot \|$ is bounded, i.e.,
\[
\sup \left\{ f(x)/\|x \| : x\in B(0,1)\setminus \{0\} \right\} <\infty.
\]

\item
$f(0)=0$.
\end{enumerate}

\paragraph{Embeddings.}
Let $C(\unitSphere)$ denote the Banach space of all continuous functions $v$ on $\unitSphere$ with respect to the norm $\|v\|_C = \sup_{x\in\unitSphere} |v(x)|.$ Let $\overline{\mathbf{C}}:=\overline{C} ([0,1],C(\unitSphere))$ be the set of all functions $f: [0,1]\to C(\unitSphere)$ such that $f$ is bounded, left continuous with respect to $\alpha\in (0,1]$, right continuous at 0, and $f$ has right limit for any $\alpha\in (0,1)$. Then we have that $\overline{\mathbf{C}}$ is a Banach space with the norm $\|f\|_{\overline{C}} = \sup_{\alpha\in [0,1]} \|f(\alpha) \|_C$.
\\
Let $\mathcal{L}:=L^2[[0,1]\times\unitSphere ; \mathbb{R}]$ be the Hilbert space of square integrable real--valued functions defined on $[0,1]\times\unitSphere$.
\\
It is known, cf. \cite{li:ogu:kre02, nat00, ram:col:gon:gil10}, that the injection $j$ defined by
\begin{equation}\label{eq:embedding}
\begin{array}{rccl}
j:  & \Fuzzy[kc]{\Banach} & \to & \overline{\mathbf{C}}\cap \mathcal{L} \\
 & \nu & \mapsto & j(\nu) = \support{\nu},
\end{array}
\end{equation}
satisfies the following properties:
\begin{enumerate}
\item
$j(r\nu^1 \Mink t\nu^2) = rj(\nu^1) + tj(\nu^2)$, $\nu^1, \nu^2\in\Fuzzy[kc]{\Banach}$ and $r,t \ge 0$.

\item
$j$ is an isometric mapping, i.e. for every $\nu^1, \nu^2\in\Fuzzy[kc]{\Banach}$,
\[
\HausDist^\infty (\nu^1, \nu^2) = \| j(\nu^1) - j(\nu^2) \|_{\overline{C}},
\qquad\textrm{and}\qquad
d_2 (\nu^1, \nu^2) = \| j(\nu^1) - j(\nu^2) \|_{\mathcal{L}}.
\]
\end{enumerate}

\paragraph{Fuzzy random variables.}
Let $(\Omega,\salgebra,\prob)$ be a complete probability space. A {\em fuzzy set--valued random variable} (\frv{}) is a function $X:\Omega \to \Fuzzy[kc]{\Banach}$, such that $X_\alpha: \omega\mapsto X(\omega)_\alpha$ are random compact convex sets for every $\alpha\in (0,1]$ (i.e. $X_\alpha$ is a $\Closed[kc]{\Banach}$--valued function measurable w.r.t.\@ $\borel[ {\Closed[kc]{\Banach}} ]$, the Borel $\sigma$--algebra on $\Closed[kc]{\Banach}$ generated by the metric $\HausDist$).
It has been proven in \cite{col:dom:lop:ral02} that this measurability definition is equivalent to the $\borel{(\Fuzzy[kc]{\Banach},d_2)}$--measurability and, it is necessary (but not sufficient) for the $\borel{(\Fuzzy[kc]{\Banach},\HausDist^\infty)}$--measurability, where $\borel{(\Fuzzy[kc]{\Banach}, D)}$ denotes the Borel $\sigma$--algebra defined on $\Fuzzy[kc]{\Banach}$ w.r.t.\@ the metric $D$.
\\
As a consequence of continuity of $\genSteiner (\cdot)$, if $X$ is a \frv{}, then $\genSteiner (X)$ is a random element in $\Banach$.
\\%
A \frv{} $X$ is \emph{integrably bounded} and we write $X\in L^1[\Omega;\Fuzzy[kc]{\Banach}]$, if $\mathbb{E} (\sup_{x\in X_{0}} \|x\| )<+\infty $.
The (Aumann) \emph{expected value} of $X\in L^1[\Omega;\Fuzzy[kc]{\Banach}]$, denoted by $\mathbb{E}[X]$, is a fuzzy set such that, for every $\alpha\in[0,1]$,
\[
(\mathbb{E}[X])_\alpha = \int_\Omega X_\alpha \dd \mathbb{P}   =  \{\mathbb{E}(f): f\in L^1[\Omega; \Banach], f \in X_\alpha \ \mathbb{P}-\textrm{a.e.}\}.
\]
It should be pointed out that, whenever $\mathbb{E} [(\sup_{x\in X_{0}} \|x\| )^2 ]< +\infty $ (we write $X\in L^2[\Omega;\Fuzzy[kc]{\Banach}]$), the expected value in the Aumann's sense is even the Frèchet expectation with respect to $d_2$, i.e.
\[
\mathbb{E}X= \underset{U\in\Fuzzy[kc]{\Banach}}{\operatorname{arg\, min}}\, \mathbb{E}(d_2(X,U)^2),
\]
see for example \cite{nat97}.
\\
In view of above measurability consideration and from embedding \eqref{eq:embedding} it follows that every \frv{} $X$ can be regarded as a random element in 
$\mathcal{L}$
, where $\support{X}(\cdot,\cdot)(\omega)= \support{X(\omega)}(\cdot,\cdot)$.
Moreover, if $X\in L^1[\Omega;\Fuzzy[kc]{\Banach}]$, for any $(x, \alpha)\in \Banach\times [0,1]$, $\support{X(\cdot)}(x,\alpha)\in L^1[\Omega;\mathbb{R}]$ and
\begin{equation}\label{eq:support_and_expectation}
\mathbb{E}[\support{X}(x,\alpha)] = \support{\mathbb{E}X}(x,\alpha).
\end{equation}

Finally, let $L^2[\Omega;\mathcal{L}] := \{ f:\Omega\to \mathcal{L} \textrm{ s.t. } [\int_\Omega \|f(\omega)\|^2_\mathcal{L} \dd \mathbb{P}]^{1/2}< +\infty \}$. It is easy to show that the map
\[
\begin{array}{rccl}
J:  & L^2[\Omega;\Fuzzy[kc]{\Banach}] & \to & L^2[\Omega;\mathcal{L}] \\
 & X & \mapsto & J(X) = j(X(\cdot)) = \support{X(\cdot)},
\end{array}
\]
is well--defined and induces an isometry in the following sense: for every $X^1, X^2\in L^2[\Omega;\Fuzzy[kc]{\Banach}] $,
\[
\Delta_2(X^1, X^2): = \mathbb{E} ( d_2 (X^1, X^2) ) = \mathbb{E} (\| J(X^1) - J(X^2) \|_{\mathcal{L}}).
\]

\section{Hukuhara set}
\label{sec:hukuhara_set}

In this section we define the \emph{Hukuhara set} associated to a \frv{} $X$, namely $\huku{X}$. We provide some properties of $\huku{X}$ most of which turn out to be useful in the next section where a decomposition theorem for fuzzy random variables is set.

Let $K$ be in $\Fuzzy[kc]{\Banach}$ such that $\genSteiner(K)=0$ and consider
\[
\detHuku{K}=\{B\in\Fuzzy[kc]{\Banach} : \genSteiner(B)=0 \textrm{ and } \exists A\in\Fuzzy[kc]{\Banach} \textrm{ s.t. } B\Mink A = K \};
\]
i.e.\ the family of those centered convex compact fuzzy sets $B$ for which the Hukuhara difference $K\HukuDiff B$ does exist.
Note that $\detHuku{K}$ is not empty, since $\indicator{0},K\in \detHuku{K}$ and $\{\lambda\odot K\}_{\lambda\in[0,1]} \subseteq \detHuku{K}$. Clearly, if $B\in\detHuku{K}$ and $A$ is the Hukuhara difference between $K$ and $B$, then $A\in\detHuku{K}$.

\begin{pro}\label{pro:detHuku_closed}
$\detHuku{K}$ is a closed subset in ($\Fuzzy[kc]{\Banach},\HausDist^\infty$).
\end{pro}
\begin{myproof}
Let $\{B_n\}\subset \detHuku{K}$ be a convergent sequence with limit $B\in\Fuzzy[kc]{\Banach}$ with respect to $\HausDist^{\infty}$, we have to prove that $B\in\detHuku{K}$. Equivalently, we have to prove that there exists $A\in \Fuzzy[kc]{\Banach}$ such that $B\Mink A = X$.
For each $n=1,2,\ldots$ there exist $A_n\in\Fuzzy[kc]{\Banach}$ such that $B_n\Mink A_n = K$.
Thus, the idea is to prove that $\{A_n\}_{n=1}^\infty$ converges, w.r.t. $\HausDist^{\infty}$, to some $A\in \Fuzzy[kc]{\Banach}$ such that $B\Mink A=X$.
To do this, let us consider the following chains of equalities
\begin{align*}
\HausDist^\infty (A_m , A_n)
&=
\| \support{A_m} - \support{A_n} \|_{\overline{C}}
\\
&=
\|(\support{A_m} + \support{B_m}) - (\support{A_n} +\support{B_n}) +\support{B_n} - \support{B_m}\|_{\overline{C}}
\\
&=
\| \support{K} - \support{K} +\support{B_n} - \support{B_m}\|_{\overline{C}}
\\
&=
\| \support{B_n} - \support{B_m}\|_{\overline{C}} = \HausDist ^{\infty} (B_n,B_m) \to 0, \quad \textrm{ for } n,m\to \infty
\end{align*}
where we use the isometry $A\mapsto \support{A}$ (first and last equalities) and the fact that $B_n,B_m$ belong to $\detHuku{K}$ (third equality).
Above limit implies that $\{ {A_n} \}_{n\in\mathbb{N}}$ is a Cauchy sequence in $(\Fuzzy[kc]{\Banach},\HausDist^{\infty})$ that is a complete metric space (e.g.~\cite[Theorem~5.1.6]{li:ogu:kre02}), and then there exists $A$ in $\Fuzzy[kc]{\Banach}$ such that $A_n\to A$.
As a consequence, $B_n\Mink A_n\to B\Mink A$ for $n\to\infty$ combined with
\[
0= \HausDist^\infty (B_n\Mink A_n, X), 
\]
guarantees that  $B\Mink A = X$ and hence $B\in\detHuku{X}$; that is the thesis.
\end{myproof}

In what follows we need the next lemma according to which a fuzzy set can be defined starting from its $\alpha$-cuts.
\begin{lem}[ (See \protect{\cite[Proposition~6.1.7, p.39]{dia:klo94}}) ]
\label{lem:charact_fuzzySet_via_alpha_cuts}
If $\{C_\alpha\}_{\alpha\in [0,1]}$ satisfies
\begin{enumerate}[(a)]
\item
$C_\alpha$ is a non empty compact convex subset of $\Banach $, for every $\alpha\in [0,1]$;

\item
$C_\beta\subseteq C_\alpha$ for $0\le \alpha\le \beta \le 1$;

\item
$C_\alpha = \bigcap_{i=1}^\infty C_{\alpha_i}$ for all sequence $\{\alpha_i\}_{i\in\mathbb{R}}$ in $[0,1]$ converging from below to $\alpha$, i.e. $\alpha_i\uparrow \alpha$ in $[0,1]$;
\end{enumerate}
then the function
\[
\nu(x)=
\left\{
\begin{array}{ll}
0,  & \textrm{if } x \not \in C_0,
\\
\sup \{ \alpha\in [0,1] : x\in C_\alpha \}, & \textrm{if } x\in C_0,
\end{array}
\right.
\]
is an element of $\Fuzzy[kc]{\Banach}$ with $\nu_\alpha = C_\alpha$ for any $\alpha\in (0,1]$ and
\[
\nu_{0} = \overline{ \bigcup_{\alpha\in (0,1]} C_\alpha } \subseteq C_0.
\]
\end{lem}

Let $X$ be a \frv{}. For the sake of simplicity and without loss of generality, let us suppose that $\genSteiner(X)=0$; otherwise one can always considered its associated centered \frv{} $\tilde{X}=X-\indicator{\genSteiner(X)}$. Next theorem defines the \emph{Hukuhara set} $\huku{X}$ associated to $X$, and provides some properties of $\huku{X}$.
\begin{pro}\label{pro:huku_set_properties}
If $B\in\Fuzzy[kc]{\Banach}$, then $E=\{B\in \detHuku{X}\}:=\{\omega\in\Omega : B\in \detHuku{X(\omega)}\}$ is measurable in $(\Omega,\salgebra)$. Moreover, if $\huku{X}=\{ B\in\Fuzzy[kc]{\Banach} : \prob ( B\in \detHuku{X} )=1 \}$, then the following statements hold.
\begin{enumerate}[{\bf (i) }]
\item\label{pro:huku_set_empty}
$\huku{X}$ is non--empty.

\item\label{pro:B_in_huku_admit_frv}
$B\in \huku{X}$ if and only if there exist a \frv{} $A$ such that $B\Mink A = X$, $\mathbb{P}$--a.s.. If $X\in L^2[\Omega; \Fuzzy[kc]{\Banach}]$, then $A$ is in $L^2[\Omega; \Fuzzy[kc]{\Banach}]$ too.

\item\label{pro:huku_set_convex}
$\huku{X}$ is a convex subset in $(\Fuzzy[kc]{\Banach},\Mink)$. As a consequence, if $B\in\huku{X}$, then $\{\lambda B\}_{\lambda\in[0,1]}\subseteq \huku{X}$.

\item\label{pro:huku_set_closed_in_HausDist}
$\huku{X}$ is a closed subset of $(\Fuzzy[kc]{\Banach},\HausDist^\infty)$.

\item\label{pro:huku_set_closed_in_dp,2}
$\huku{X}$ is a closed subset of $(\Fuzzy[kc]{\Banach},d_{2})$.
\end{enumerate}
\end{pro}
\begin{myproof}
Using the definition of $\detHuku{X(\omega)}$ and the characterization of element in $\Fuzzy[kc]{\Banach}$ via the support functions, we get the following chains of equalities.
\begin{align*}
E & = \{\omega\in\Omega : \genSteiner(B)=0 \textrm{ and } \exists A_\omega\in\Fuzzy[kc]{\Banach}, \ B\Mink A_\omega =X(\omega)\}
\\
& = \{ \omega\in\Omega : \genSteiner(B)=0 \} \cap \{\omega\in\Omega : \exists A_\omega\in\Fuzzy[kc]{\Banach}, \textrm{ s.t. }\ \support{B}+ \support{A_\omega} =\support{X(\omega)}\}.
\end{align*}
Since $B$ is a deterministic fuzzy set, $E_0=\{ \omega\in\Omega : \genSteiner(B)=0 \}$ is either the empty set or the whole $\Omega$; hence $E_0$ is measurable. On the other hand, $A_\omega$ in $\Fuzzy[kc]{\Banach}$ satisfies $B\Mink A_\omega=X(\omega)$ if and only if $\support{B}+\support{A_\omega} = \support{X(\omega) }$ or, equivalently, if and only if $\support{X(\omega) }-\support{B}$ is the support function of some element in $\Fuzzy[kc]{\Banach}$. Thus, because of Properties~1--6 we have that
\begin{align*}
E & = E_0 \cap \{\omega\in\Omega : \ f_{\omega} \textrm{ satisfies Properties 1--6}\}
\\
& = E_0 \cap E_1\cap \ldots \cap E_6,
\end{align*}
where $E_i=\{\omega\in\Omega: f_\omega \textrm{ satisfies Property }i\}$ for $i=1,\ldots,6$. If $E_1,\ldots,E_6$ are measurable events, then $E$ is measurable too. To show this note that each $E_i$ ($i=1,\ldots, 6$) can be written as $E_i=\{\omega: g_i(\omega)\le 0 \}$ where
\begin{align*}
g_1 &= \sup \{ | \limsup_{y\to x} f_\omega(y)-f_\omega(x) | : {x\in B(0,1)} \},
\\
g_2 &= \sup \{ \lambda f_\omega(x)-f_\omega(\lambda x)  :\lambda\in (0,1], x\in B(0,1) \},
\\
g_3 &= \sup \Big\{
\|x\| f_\omega\left(\lambda \frac{x}{\|x\|}\right) - \|x_1\| f_\omega\left(\lambda \frac{x_1}{\|x_1\|}\right) -
\|x_2\| f_\omega\left(\lambda \frac{x_2}{\|x_2\|}\right)
\\
& \quad\qquad\qquad :\lambda\in (0,1], x,x_1,x_2\in \Banach\setminus\{0\}, \textrm{ with } x=x_1+x_2
\Big\},
\\
g_4 &= - \sup \{ f_\omega(x) + f_\omega(-x) : x\in B(0,1)\},
\\
g_5 &= \sup \Big\{ \frac{|f_\omega(x)|}{\|x\|} : x\in B(0,1)\setminus\{0\}\Big\},
\\
g_6 &= |f_\omega(0)|.
\end{align*}
Clearly $\omega\mapsto g_i(\omega)$ are measurable maps and hence $E$ is a measurable event in the $\sigma$--algebra $\salgebra$.
\\
\noindent {\sl \bf ITEM (\ref{pro:huku_set_empty}).}
Surely $\indicator{0}$ belongs to $\huku{X}$, hence $\huku{X}$ is not empty.

\noindent {\sl \bf ITEM (\ref{pro:B_in_huku_admit_frv}).}
The sufficiency is trivial, let us prove the necessity. Let $E^c=\Omega\setminus E = \{\omega\in\Omega : B\not\in\detHuku{X(\omega)} \}$,
by hypothesis $\mathbb{P}(E)=1$ and $\mathbb{P}(E^c)=0$. For every $\omega\in\Omega\cap E$, there exists $A_\omega\in \Fuzzy[kc]{\Banach}$ such that $B\Mink A_\omega = X(\omega)$. Let us consider the map
\begin{equation}\label{eq:A_omega_defines_an_FRV}
\begin{array}{rccl}
A:  & \Omega & \to & \Fuzzy[kc]{\Banach} \\
 & \omega & \mapsto & A(\omega)=
\left\{
\begin{array}{ll}
A_\omega, & \omega\in\Omega\cap E,
\\
\indicator{0}, & \omega\in E^c.
\end{array}
\right.
\end{array}
\end{equation}
Since $\support{A}=\support{X}-\support{B}$ $\mathbb{P}$--almost surely, $\support{A}$ is measurable. Hence, the map $A$ defined above, is the \frv{} we are looking for.
\\%
Moreover, let $X\in L^2[\Omega;\Fuzzy[kc]{\Banach}]$, then $\support{X}$ and hence $\support{A}=\support{X}-\support{B}$ belong to $L^2[\Omega;\mathcal{L}]$.

\noindent {\sl \bf ITEM (\ref{pro:huku_set_convex}).}
Consider $B_1,B_2\in\huku{X}$. From above part we know that there exist two \frv{} $A_1,A_2$ with values in $\Fuzzy[kc]{\Banach}$ such that $\mathbb{P}$--a.s. $B_1\Mink A_1=X$ and $B_2\Mink A_2=X$. For any $\lambda\in [0,1]$, the following hold
\[
\lambda (B_1\Mink A_1)=\lambda X, \qquad (1-\lambda) (B_2\Mink A_2) =(1-\lambda) X, \qquad \mathbb{P}-a.s.
\]
from which we get
\[
\lambda B_1\Mink  (1-\lambda) B_2 \Mink  A = X, \qquad \mathbb{P}-a.s.
\]
with $A= \lambda A_1\Mink  (1-\lambda) A_2$ $\mathbb{P}$--a.s.\@. Hence $\lambda B_1\Mink  (1-\lambda) B_2 \in \huku{X}$.
\\
To prove the last part consider $B\in\huku{X}$, then $\lambda B = \lambda B \Mink (1-\lambda) \indicator{0} \in\huku{X}$.

\noindent {\sl \bf ITEM (\ref{pro:huku_set_closed_in_HausDist}).}
Consider a sequence $\{B_n\}_{n=1}^\infty\subset \huku{X}$ converging to $B\in\Fuzzy[kc]{\Banach}$ in $(\Fuzzy[kc]{\Banach},\HausDist^\infty)$, i.e.
\[
\HausDist ^{\infty} (B,B_n) \to 0,  \textrm{ as } n\to\infty.
\]
We have to prove that $B\in \huku{X}$. For any $n\in\mathbb{N}$, let $E_n=\{\omega\in\Omega : B_n \in\huku{X}\}$ and $A_n$ a \frv{} as in (\ref{pro:B_in_huku_admit_frv}). Then for every $\omega\in\Omega\cap E_n$, $B_n\Mink A_{n}(\omega) = X(\omega)$ and
\[
\HausDist^\infty (A_{m}(\omega) , A_{n}(\omega)) = \HausDist ^{\infty} (B_m,B_n) \to 0, \textrm{ as }n\to\infty.
\]
Thus, the completeness of $(\Fuzzy[kc]{\Banach},\HausDist^\infty)$ guarantees that, for every $\omega\in\Omega\setminus \bigcup_n (E_n)^c = \Omega\cap \bigcap_n E_n$, $\{A_{n}(\omega)\}_{n\in\mathbb{N}}$ converges w.r.t. $\HausDist^\infty$ to some $A_\omega\in\Fuzzy[kc]{\Banach}$. Further, for every $\omega\in\Omega\cap \bigcap_n E_n$ and $n\in\mathbb{N}$ the following inequalities hold
\begin{align*}
0\le \HausDist^\infty (X(\omega),B\Mink A_\omega) & \le \HausDist^\infty (X(\omega),B_n \Mink A_{n}(\omega)) + \HausDist^\infty (B_n \Mink A_{n}(\omega),B\Mink A_\omega)
\\
& \le 0 + \HausDist^\infty (B_n,B) + \HausDist^\infty (A_{n}(\omega),A_\omega) \to 0
\end{align*}
where, for the first addend, we use the fact that $X(\omega)=B_n \Mink A_{n}(\omega)$. Then $X = B \Mink A$ {} $\mathbb{P}$--a.s., and $A$ is the \frv{} defined by Equation~\eqref{eq:A_omega_defines_an_FRV}. Thus we have the thesis; the limit of the convergent sequence $\{B_n\}\subseteq \huku{X}$ belongs to $\huku{X}$ too.

\noindent {\sl \bf ITEM (\ref{pro:huku_set_closed_in_dp,2}).}
Let us consider a sequence $\{B_n\}_{n=1}^\infty\subset \huku{X}$ converging to $B\in\Fuzzy[kc]{\Banach}$ in $(\Fuzzy[kc]{\Banach},d_{2})$, i.e.
\[
d_{2} (B,B_n) \to 0,  \textrm{ as } n\to\infty.
\]
We have to prove that $B\in \huku{X}$. In this case, $(\Fuzzy[kc]{\Banach},d_{2})$ is not complete and, hence, we can not repeat all arguments in (\ref{pro:huku_set_closed_in_HausDist}). In particular, for any $n\in\mathbb{N}$ and for every $\omega\in\Omega\cap E_n=\{\omega\in\Omega : B_n\in\detHuku{X(\omega)} \}$, there exist $A_{n}(\omega)$ such that $B_n\Mink A_{n}(\omega) = X(\omega)$ and, using analogous arguments of those in Proposition~\ref{pro:detHuku_closed},
\[
d_{2} (A_{m}(\omega) , A_{n}(\omega)) = \left(\int_0^1 \int_{\unitSphere} | \support{A_{m}(\omega)}(\alpha, u) - \support{A_{n}(\omega)}(\alpha, u) |^2 \dd u \dd \alpha \right)^{\frac{1}{2}} = d_{2} (B_m,B_n) \to 0,
\]
as $n\to\infty$ and where $\dd\alpha$ and $\dd u$ denote the Lebesgue measure on $[0,1]$ and the normalized Lebesgue measure on $\unitSphere$ respectively.
Thus, for every $\omega\in\Omega\cap \bigcap_n E_n$, $\{\support{A_{n}(\omega)}\}_{n\in\mathbb{N}}$ is a Cauchy sequence in the Hilbert space $\mathcal{L}$ ($= L^2[[0,1]\times \unitSphere ; \mathbb{R}]$) and it admits limit in $\mathcal{L}$, namely $f_\omega$.
Since
\[
\|\support{A_{n}(\omega)} - (\support{X(\omega)}- \support{B}) \|_{\mathcal{L}} =
\|(\support{A_{n}(\omega)} - \support{X(\omega)}) + \support{B} \|_{\mathcal{L}} =
\|\support{B} -\support{B_{n}} \|_{\mathcal{L}} \to 0,
\]
necessarily we have
\[
\support{A_{n}(\omega)} \stackrel{L^2}{\to} f_\omega = \support{X(\omega)} - \support{B},\qquad \forall \omega\in\Omega\cap \bigcap_n E_n.
\]
Note that, $f_\omega$ is not necessarily the support function of some element in $\Fuzzy[kc]{\Banach}$.
In other words, for every $\omega\in\Omega\cap \bigcap_n E_n$, $\{A_{n}(\omega)\}_{n\in\mathbb{N}}$ is a Cauchy sequence in the non--complete space $(\Fuzzy[kc]{\Banach},d_{2})$, but under the embedding $j$, Equation~\eqref{eq:embedding}, we have that the sequence $\{j(A_{n}(\omega))\}_{n\in\mathbb{N}} =\{\support{A_{n}(\omega)}\}_{n\in\mathbb{N}}$ is a Cauchy sequence that admits limit in the Hilbert space $\mathcal{L}$. But, in general, this limit is not the image under $j$ of some element of $\Fuzzy[kc]{\Banach}$.
We claim that, for every $\omega\in\Omega\cap \bigcap_n E_n$, there exists $A_\omega\in\Fuzzy[kc]{\Banach}$ such that $\support{A_\omega}=f_\omega = \support{X(\omega)} - \support{B}$. This allows us to deduce the thesis because, defining the \frv{} $A$ as in Equation~\eqref{eq:A_omega_defines_an_FRV}, we have that $B\Mink A = X$ holds $\mathbb{P}$--a.s.\@.
\\
In fact, let us consider the family $\{C_\alpha\}_{\alpha\in [0,1]}$ of subsets of ${\Banach}$ defined by
\[
C_\alpha  = \{ y\in \Banach  : \langle y,u\rangle \le f_\omega(\alpha,u), \forall u\in\unitSphere \},\quad \alpha\in[0,1].
\]
In what follows, let $\omega\in\Omega\cap\bigcap_n E_n$, we prove that the family $\{C_\alpha\}_{\alpha\in [0,1]}$ satisfies (a), (b), (c) from Lemma~\ref{lem:charact_fuzzySet_via_alpha_cuts}, and it defines uniquely a fuzzy set $\nu$ whose support function is, clearly, $f_\omega$. Thus the fuzzy set $\nu$ defined in Lemma~\ref{lem:charact_fuzzySet_via_alpha_cuts} is just the $A(\omega)$ in $\Fuzzy[kc]{\Banach}$ we are looking for.
\\
(a). Let $\alpha \in [0,1]$.
\\
\emph{$C_\alpha$ is non--empty}: since $B_\alpha\subseteq (X(\omega))_\alpha$, then for every $u\in\unitSphere$
\begin{equation}\label{eq:f_omega=hX-hB}
f_\omega(\alpha,u) = \support{X(\omega)}(\alpha,u) - \support{B} (\alpha,u) \ge 0 = \langle 0, u\rangle,
\end{equation}
i.e. $0\in C_\alpha$.
\\
\emph{$C_\alpha$ is convex}: let $\lambda\in[0,1]$ and $y_1,y_2\in C_\alpha$, for every $u\in\unitSphere$
\[
\langle \lambda y_1 + (1-\lambda) y_2 , u \rangle \le \lambda f_\omega (\alpha,u) + (1- \lambda ) f_\omega (\alpha,u) =
f_\omega (\alpha,u)
\]
i.e. $\lambda y_1 + (1-\lambda) y_2 \in C_\alpha$.
\\
\emph{$C_\alpha$ is compact}: we have to prove that it is a bounded closed subset of $\Banach $. Note that $\{0\}\subseteq B_\alpha \subseteq (X(\omega))_\alpha$, then $\support{X(\omega)}(\alpha,u)\ge \support{B}(\alpha,u)\ge 0$ for each $u\in\unitSphere$ and $\support{X(\omega)}(\alpha,u)\ge \support{X(\omega)}(\alpha,u) - \support{B}(\alpha,u) = f_\omega (\alpha, u)$. This implies that $\langle y,u \rangle $ is bounded for every $u\in\unitSphere$ and hence that $C_\alpha \subseteq \Banach $ is bounded. On the other hand, let $\{y_n\}\subset C_\alpha$ be convergent to $y\in\Banach $, then, for every $n\in\mathbb{N}$ and $u\in\unitSphere$,
\[
\langle y_n,u \rangle \le f_\omega(\alpha,u),
\]
and passing to the limit we obtain the same inequality for $y$ and for every $u\in\unitSphere$; i.e. $y\in C_\alpha$. This fact allows us to conclude that $C_\alpha$ is closed and hence compact.
\\
(b). Let $0\le \alpha\le \beta \le 1$.
Note that, for every $n\in\mathbb{N}$ and $u\in\unitSphere$, $\support{A_{n}(\omega)}(\beta,u) \le \support{A_{n}(\omega)}(\alpha,u)$. Let $n\to \infty$, then $f_\omega(\beta,u) \le f_\omega(\alpha,u)$ for every $u\in\unitSphere$; i.e., for every $u\in\unitSphere$ and $n\in\mathbb{N}$, $\support{A_{n}(\omega)}$ and $f_\omega$ are non--increasing functions with respect to $\alpha$.
Now, let us consider $y\in C_\beta$, then for every $u\in\unitSphere$, $\langle y, u \rangle \le f_\omega (\beta,u)\le f_\omega (\alpha, u)$; i.e. $y\in C_\alpha$ and $C_\beta\subseteq C_\alpha$.
\\
(c). Let $\{\alpha_i\}_{i\in\mathbb{N}}\subset [0,1]$ such that $\alpha_i\uparrow \alpha$ as $i$ tends to infinity, that is $\alpha_i\le \alpha_{i+1}$ and $\alpha_i\to \alpha$ as $i\to\infty$. Because of $\alpha_i\le \alpha$ and (b), we have $C_{\alpha}\subseteq C_{\alpha_i}$ and $C_\alpha \subseteq \bigcap_{i\in\mathbb{N}} C_{\alpha_i}$. It remains to show the opposite inclusion. To do this let $y\in \bigcap_{i\in\mathbb{N}} C_{\alpha_i}$, i.e. $y\in C_{\alpha_i}$ for all $i\in\mathbb{N}$ or, equivalently,
\begin{equation}\label{eq:y_in_intersection_of_C_alpha}
\langle y,u\rangle\le f_\omega(\alpha_i,u),\qquad \textrm{for every } i\in\mathbb{N}, u\in\unitSphere.
\end{equation}
Note that, for every $u\in\unitSphere$, $f_\omega(\cdot,u)$ is left--continuous with respect to $\alpha$ because it is the difference of two left--continuous functions (cf. Equation~\eqref{eq:f_omega=hX-hB}). Hence, for the arbitrariness of $i$ in \eqref{eq:y_in_intersection_of_C_alpha}, as $i$ tends to infinity we get $\langle y,u\rangle\le f_\omega(\alpha,u)$; i.e. $y\in C_\alpha$.
\end{myproof}

\section{Hukuhara decomposition}
\label{sec:hukuhara_decomposition}

Let us recall again the well known decomposition \eqref{eq:gauss_decomposition} for Gaussian \frv{} $X$
\[
X=\mathbb{E} X \Mink \indicator{\xi },
\]
where $\mathbb{E} X$ is the Aumann expectation of $X$, and $\xi$ is a Gaussian random element in $\Banach$ with $\mathbb{E} \xi=0$. Equation~\eqref{eq:gauss_decomposition} implies a \emph{randomness defuzzification} for \frv{} $X$ that is equal to its expected value $\mathbb{E} X$ (a deterministic fuzzy set) up to a random Gaussian translation $\xi$.
In \cite{bon11}, the author showed another case of defuzzification of randomness: a Brownian fuzzy set--valued process is reduced to be a Brownian process in $\Banach$.
In both cases, the randomness initially defined on $\Fuzzy[kc]{\Banach}$ can be simply defined on $\Banach$.
\\
Now, our question becomes the following one. \emph{Under what conditions can we establish that a defuzzification of randomness occurs for fuzzy set--valued random process?} In other words: \emph{Can a fuzzy process, whose randomness is given only by vectors, be characterized in some way?}
In this section we propose a positive answer to the above question. We focus mainly on a decomposition theorem for \frv{}. In fact, in Theorem~\ref{teo:X=C+Y}, we prove that any \frv{} $X$ can be decomposed as the sum of a deterministic convex fuzzy set $\hukuMin{X}$ and a \frv{} $Y$ (that contains the whole randomness) in a unique way. This decomposition allows us to characterize, by means of the Aumann expected value, the \frv{} that is a random translation of a deterministic fuzzy set.
\begin{defi}
A \frv{} $X$ is a \emph{translation} if there exists $M\in\Fuzzy[kc]{\Banach}$ with $\genSteiner (M)=0$ such that
\[
X(\omega)=M\Mink\indicator{\genSteiner(X)}.
\]
\end{defi}
Roughly speaking, the randomness of a translation depends only on the specific location in the underling space $\Banach$ while it does not depend on its fuzzy shape. Note that, accordingly to \eqref{eq:gauss_decomposition}, every Gaussian \frv{} $X$ is a \frv{} translation with $M\Mink \indicator{\mathbb{E}(\genSteiner (X) )}=\mathbb{E}X$. Another sufficient condition for $X$ to be a \frv{} translation is given by Proposition~\ref{pro:degenere_implies_translation}, while a necessary and sufficient condition is state in Theorem~\ref{teo:transl_iff_C=EX}.
\begin{pro}\label{pro:degenere_implies_translation}
Let $X$ be a \frv{} such that $\mathbb{E} X=\indicator{c}$ where $c\in\Banach$. Then $X=\indicator{\xi}$ $\mathbb{P}$--a.s.\@ for some random element $\xi$ in $\Banach$. (Clearly $X$ is a \frv{} translation.)
\end{pro}
\begin{myproof}
Thesis can be obtained using similar arguments in \cite[Theorem~8]{bon11}, or, whenever $X\in L^2[\Omega;\Fuzzy[kc]{\Banach}]$, as corollary of the Theorem~\ref{teo:X=C+Y} and Theorem~\ref{teo:transl_iff_C=EX}.
\end{myproof}
Clearly, the vice versa of Proposition~\ref{pro:degenere_implies_translation} does not hold, for example in the case of Gaussian \frv{}.
In order to characterize translation \frv{}, we need the following decomposition theorem.
\begin{teo}\label{teo:X=C+Y}
Let $X\in L^2[\Omega;\Fuzzy[kc]{\Banach}]$ with $\genSteiner(X)=0$. Thus there exists $\hukuMin{X}\in\Fuzzy[kc]{\Banach}$ with $\genSteiner (\hukuMin{X})=0$ and $Y\in L^2[\Omega;\Fuzzy[kc]{\Banach}]$ such that $X$ decomposes according to
\begin{equation}\label{eq:X=C+Y}
X(\omega)=\hukuMin{X}\Mink Y(\omega),
\end{equation}
for $\mathbb{P}$--almost all $\omega\in\Omega$. In particular, $\hukuMin{X}$ is the unique element in $\Fuzzy[kc]{\Banach}$ that satisfies \eqref{eq:X=C+Y} and minimizes $\mathbb{E}[ (d_2(X,C))^2 ]$; i.e., there exists a unique $\hukuMin{X}\in\huku{X}$ such that
\begin{equation}\label{eq:C_argument_that_minimizes}
\hukuMin{X}:=\mathop{arg\, min}_{B\in\huku{X}} \mathbb{E}[ (d_2(X,B))^2 ].
\end{equation}
Hence $Y$ is the unique (except on a $\mathbb{P}$--negligible set) \frv{} such that its support function is given by $\support{Y} = \support{X}-\support{\hukuMin{X}}$.
Moreover, $\hukuMin{X}$ is a maximal element in $\huku{X}$ with respect to the level--wise set inclusion; that is, if $C\in\huku{X}$ with $(\hukuMin{X})_\alpha \subseteq C_\alpha$ for any $\alpha\in [0,1]$, then  $\hukuMin{X}=C$.
\end{teo}
\begin{myproof}
Since $\huku{X}$ collects all the element of $\Fuzzy[kc]{\Banach}$ for which \eqref{eq:X=C+Y} holds, we have to prove that there exists a unique element in $\huku{X}$ that minimizes the map $B\in\huku{X}\to \mathbb{E}[ (d_2(X,B))^2 ]$.
\\
At first note that $\huku{X}$ can be seen as a subset of $L^2[\Omega; \Fuzzy[kc]{\Banach}]$; in fact, for each $B\in\huku{X}$ the constant map $\omega\mapsto B$ is an element of $L^2[\Omega; \Fuzzy[kc]{\Banach}]$ since
\[
\mathbb{E}[ (\sup_{b\in B_0} \|b\|)^2 ] = (\sup_{b\in B_0} \|b\|)^2 <+\infty.
\]
Moreover, $\huku{X}$ is closed in $L^2[\Omega;\Fuzzy[kc]{\Banach}]$ as a consequence of
\[
\mathbb{E}[ (d_2(A,B) )^2] = (d_2(A,B) )^2,
\]
for any couples $A,B\in\Fuzzy[kc]{\Banach}$, and thanks to the fact that $\huku{X}$ is closed in $(\Fuzzy[kc]{\Banach},d_2)$, see Proposition~\ref{pro:huku_set_properties}.
\\
Thus the minimization problem is equivalent to prove that there exists a unique projection of $X$ onto $\huku{X}$ that is a closed convex subset of $L^2[\Omega;\Fuzzy[kc]{\Banach}]$ endowed with the metric $\Delta_2$. Since $L^2[\Omega;\Fuzzy[kc]{\Banach}]$ embeds isometrically in the Hilbert space $L^2[\Omega;\mathcal{L}]$ through map $J$ (see the Introduction), there exists a unique element $\hukuMin{X}\in\huku{X}$ that realizes the required minimum \eqref{eq:C_argument_that_minimizes}.
\\
As a consequence of $\hukuMin{X}\in\huku{X}$ and of \eqref{pro:B_in_huku_admit_frv} in Proposition~\ref{pro:huku_set_properties}, the \frv{} $Y$ is defined through its support function $\support{Y}=\support{X}-\support{\hukuMin{X}}$.
\\
Finally, let $C$ be as in the thesis; thus inclusions $(\hukuMin{X})_\alpha\subseteq C_\alpha \subseteq X_\alpha$ imply $\support{X}-\support{C} \le \support{X} - \support{\hukuMin{X}}$. Then, by definition of $\hukuMin{X}$ and $d_2$, necessarily $C=\hukuMin{X}$ holds.
\end{myproof}
The chosen notation wants to recall the line of the proof; $\hukuMin{X}$ is obtained through a projection theorem of the given \frv{} $X$ on its Hukuhara set $\huku{X}$.
Further, we want to stress out that the suffix $X$ does not mean that $\hukuMin{X}$ is random; in fact, it does not depend on $\omega$ but rather it is a deterministic element of $\Fuzzy[kc]{\Banach}$ (that is a constant element in $L^2[\Omega;\Fuzzy[kc]{\Banach}]$) that depends on the whole map $\omega\mapsto X(\omega)$.

The following theorems provide necessary and sufficient condition for a \frv{} to be a translation.
\begin{teo}\label{teo:necessary_cond_for_translation}
Let $X$ be a \frv{} translation, and $\widetilde{X} = X \Mink \indicator{-\genSteiner (X)}$. Then
\begin{equation}\label{eq:X=H_X+steiner(X)}
X=\hukuMin{\widetilde{X}}\Mink \indicator{\genSteiner (X)},\qquad \mathbb{P}-\textrm{a.s.}
\end{equation}
\end{teo}
\begin{myproof}
By hypothesis $X=M\Mink \indicator{\genSteiner (X)}$ for some $M\in\Fuzzy[kc]{\Banach}$ with $\genSteiner (M)=0$. Clearly, $\widetilde{X}=X\Mink \indicator{-\genSteiner(X)} = M$ and $\genSteiner (\widetilde{X})=0$. Thus, by Theorem~\ref{teo:X=C+Y} applied to $\widetilde{X}$, we have $M\in\huku{\widetilde{X}}$ and $\mathbb{E} [ ( d_2(M,\widetilde{X})^2 ) ]=0$; that is, $M=\hukuMin{\widetilde{X}}$.
\end{myproof}
\begin{teo}\label{teo:transl_iff_C=EX}
Let $X\in L^2[\Omega;\Fuzzy[kc]{\Banach}]$. $X$ is a \frv{} translation if and only if $\hukuMin{\widetilde{X}}$ satisfies
\begin{equation}\label{eq:transl_iff_C=EX}
\mathbb{E}X = \hukuMin{\widetilde{X}} \Mink \indicator{\mathbb{E}(\genSteiner(X) )}
\end{equation}
with $\mathbb{E}X$ being the Aumann expectation; in other words, $\hukuMin{\widetilde{X}}$ is $\mathbb{E}X$ up to a translation.
\end{teo}
\begin{myproof}
For the \lq\lq only if\rq\rq{} part, in order to obtain Equation~\eqref{eq:transl_iff_C=EX}, it is sufficient to compute the expectation in Equation~\eqref{eq:X=H_X+steiner(X)}.
\\
Consider the \lq\lq if\rq\rq{} part. For the sake of simplicity, let us assume that $\genSteiner (X)=0$, a straightforward argument extends the result in the more general case of a \frv{} with non--void $\genSteiner (X)$.
Then, in term of support functions, Equation~\eqref{eq:X=C+Y} becomes
\[
\support{X}=\support{\hukuMin{X}} + \support{Y} = \support{\mathbb{E}X} + \support{Y}, \qquad \mathbb{P}-a.s.
\]
where we use the fact that $\hukuMin{X}=\mathbb{E}X$. Computing expectation of both sides and using \eqref{eq:support_and_expectation}, we get $\support{\mathbb{E}Y}=0$. Hence $Y=\indicator{ \xi  }$ a.s. for some random element $\xi$ in $\Banach$ (cf. \cite{bon11}).
\end{myproof}
\begin{oss}
Whenever $X\in L^2[\Omega;\Fuzzy[kc]{\Banach}]$, in view of Theorem~\ref{teo:X=C+Y} and Theorem~\ref{teo:transl_iff_C=EX}, we get a proof of Proposition~\ref{pro:degenere_implies_translation}. In fact, suppose that
$\mathbb{E} X = \indicator{c}$ for some $c\in\Banach$, and compute expectation of both sides in Equation~\eqref{eq:X=C+Y}
\[
\indicator{c} = \mathbb{E}X= \hukuMin{X} \Mink \mathbb{E} Y.
\]
Hence, for any $\alpha\in[0,1]$, $(\hukuMin{X})_\alpha$ is a subset of $\{c\}$ up to a translation, that is $(\hukuMin{X})_\alpha$ is a singleton as well as $(\mathbb{E}Y)_\alpha$. Then $\hukuMin{X}=\indicator{c'}$ for some $c'\in\Banach$, i.e. $\hukuMin{X}$ is equal to $\mathbb{E}X$ up to a translation and, by Theorem~\ref{teo:transl_iff_C=EX}, $X$ is a \frv{} translation that implies $Y=\indicator{\xi}$ for some random element in $\Banach$. Finally, Equation~\eqref{eq:X=C+Y} becomes
\[
X= \hukuMin{X} \Mink Y= \indicator{c'}\Mink \indicator{\xi}= \indicator{\xi'},
\]
that is the thesis of Proposition~\ref{pro:degenere_implies_translation}.
\end{oss}
Moreover, the following results hold.
\begin{cor}
Let $X\in L^2[\Omega;\Fuzzy[kc]{\Banach}]$ with $\genSteiner (X)=0$ and $\mathbb{E} X = \hukuMin{X}$. Thus $X$ is almost surely deterministic and equal to $\hukuMin{X}$.
\end{cor}
\begin{cor}\label{cor:X'=X+D_then_C'=C+D}
Let $X\in L^2[\Omega;\Fuzzy[kc]{\Banach}]$, $D\in\Fuzzy[kc]{\Banach}$ and $X' = X\Mink D$ with $\genSteiner (X) =\genSteiner (D) = 0$ (hence $\genSteiner (X')=0$ too). Then $\hukuMin{X'}=\hukuMin{X}\Mink D$.
\end{cor}
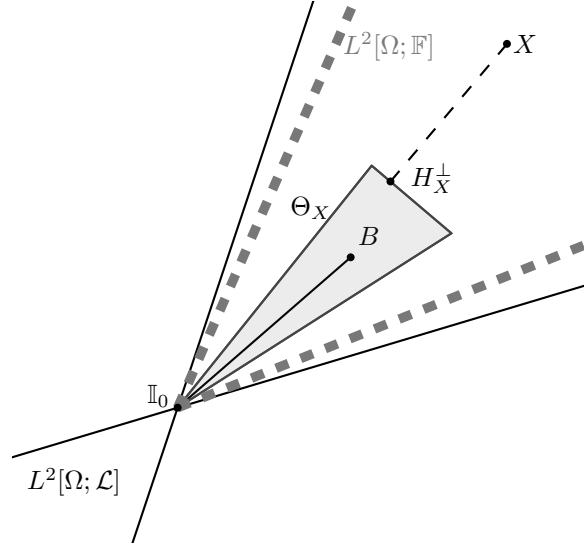
\begin{figure}[!htb]
\begin{center}
\newrgbcolor{uuuuuu}{0.27 0.27 0.27}
\newrgbcolor{uququq}{0.25 0.25 0.25}
\newrgbcolor{xxxxxx}{0.47 0.47 0.47}
\psset{xunit=1.0cm,yunit=1.0cm,algebraic=true,dotstyle=o,dotsize=3pt 0,linewidth=0.8pt,arrowsize=3pt 2,arrowinset=0.25}
\begin{pspicture*}(-2.2,-1.8)(5.4,5.4)
\rput[tl](2.4 , 2.4){$B$}
\psdots[dotstyle=*](2.3, 2)
\psplot{0}{2.3}{2/2.3*x}
\pspolygon[linecolor=uuuuuu,fillcolor=uuuuuu,fillstyle=solid,opacity=0.1](0,0)(3.64,2.32)(2.58,3.22)
\psplot{-2.7}{5.92}{(-0--0.3*x)/1}
\psplot{-2.7}{5.92}{(-0--3*x)/1}
\psline[linecolor=uuuuuu](0,0)(3.64,2.32)
\psline[linecolor=uuuuuu](3.64,2.32)(2.58,3.22)
\psline[linecolor=uuuuuu](2.58,3.22)(0,0)
\rput[tl](1.5,2.8){$\huku{X}$}
\psline[linestyle=dashed,dash=5pt 5pt](4.38,4.84)(2.83,3.01)
\rput[tl](3.1,3.25){$\hukuMin{X}$}
\rput[tl](2.2,5){${\xxxxxx{ L^2 [\Omega; \Fuzzy[kc]{\Banach}] }}$}
\psplot[linewidth=3.6pt,linestyle=dashed,dash=5pt 5pt,linecolor=xxxxxx]{0}{8.34}{(-0--2.3*x)/1}
\psplot[linewidth=3.6pt,linestyle=dashed,dash=5pt 5pt,linecolor=xxxxxx]{0}{8.34}{(-0--0.4*x)/1}
\rput[tl](-2,-0.8){${L^2 [\Omega;\mathcal{L}] }$}
\psdots[dotstyle=*](0,0)
\rput[tl](-0.4,0.3){$\indicator{0}$}
\psdots[dotstyle=*](4.38,4.84)
\rput[tl](4.46,4.96){$X$}
\psdots[dotstyle=*](2.83,3.01)
\end{pspicture*}
\caption{A qualitative graphical interpretation of some results of Section~\ref{sec:hukuhara_set} and Section~\ref{sec:hukuhara_decomposition}. In particular, $\huku{X}$ is represented as a closed convex subset of $\Fuzzy[kc]{\Banach}$ containing the origin and such that, for any $B\in\huku{X}$ and $\lambda\in [0,1]$, $\lambda B\in\huku{X}$. Hence, $\hukuMin{X}$ is the projection of $X$ on $\huku{X}$, as a subset of $L^2[\Omega;\Fuzzy[kc]{\Banach}]$, with respect to the metric $\mathbb{E}[d_2(\cdot,\cdot)^2]$, this also guarantees the uniqueness of $\hukuMin{X}$ since the cone $L^2[\Omega;\Fuzzy[kc]{\Banach}]$ is embeddable in the Hilbert space $L^2 [\Omega;\mathcal{L}]$ through the isometry $X\mapsto j(X)$. Finally the following inclusions or embeddings are qualitatively represented:
$\huku{X}\subseteq \Fuzzy[kc]{\Banach} \hookrightarrow L^2[\Omega; \Fuzzy[kc]{\Banach}] \hookrightarrow L^2 [\Omega;\mathcal{L}]$.}
\end{center}
\end{figure}
Remark~\ref{oss:example_X_not_translation} shows an example of an $X$ in $L^2[\Omega;\Fuzzy[kc]{\Banach}]$ with $\genSteiner (X)=0$ for which $\mathbb{E}(X)\neq \hukuMin{X}$ and for which $\hukuMin{X}$ is not necessarily $\indicator{0}$; i.e., in terms of Theorem~\ref{teo:transl_iff_C=EX}, $X$ is not a translation but its deterministic part $\hukuMin{X}$ in the decomposition \eqref{eq:X=C+Y} is not just reduced to the origin.
\begin{oss}\label{oss:example_X_not_translation}
Let $\Banach = \mathbb{R}$, $(\Omega = [0,1],\borel[{[0,1]}],\mathbb{P})$ where $\borel[{[0,1]}]$ denotes the Borel $\sigma$--algebra on $[0,1]$ w.r.t. the euclidean metric and $\mathbb{P}=\misura$ is the Lebesgue measure. Let $X$ be the \frv{} defined by $X:=\indicator{[\omega,\omega]}$, for any $\omega\in [0,1]$. Clearly $X\in L^2[\Omega;\Fuzzy[kc]{\Banach}]$ and $\genSteiner (X)=0$. Moreover,
\[
f_m(\omega):=\min X_1(\omega)=-\omega \quad \textrm{ and } \quad f_M(\omega):=\max X_1(\omega)=\omega
\]
are integrable selections of the 1--level \racs{} $X_1$. Obviously, any other integrable selection $f$ of $X_1$ satisfies
\[
f_m(\omega) \le f(\omega) \le f_M(\omega), \quad \textrm{ for each }\omega\in[0,1].
\]
Then
\[
-\frac{1}{2}=\mathbb{E}f_m \le \mathbb{E}f \le \mathbb{E}f_M=\frac{1}{2},
\]
and, by the convexity of Aumann expectation and because $X_1=X_\alpha$ for any $\alpha\in[0,1]$, $\mathbb{E}X_1=[-\frac{1}{2},\frac{1}{2}] = \mathbb{E}X_\alpha$, that is $\mathbb{E}X=\indicator{[-\frac{1}{2},\frac{1}{2}]}$.
\\
We prove that $\mathbb{E}X\not \in \huku{X}$ and hence, by Theorem~\ref{eq:transl_iff_C=EX}, $X$ is not a \frv{} translation. In fact, note that
\[
X\HukuDiff \mathbb{E}X = \indicator{[-\omega, \omega]} \HukuDiff \indicator{[-\frac{1}{2} , \frac{1}{2}]} =
\left\{
\begin{array}{ll}
\indicator{[-\omega + \frac{1}{2}, \omega -\frac{1}{2}]}, & \omega > \frac{1}{2},
\\
\indicator{0}, & \omega = \frac{1}{2},
\\
\textrm{it does not exist}, & \omega < \frac{1}{2},
\end{array}
\right.
\]
implies
\[
\mathbb{P} (\mathbb{E} X\in \detHuku{X} ) = \mathbb{P} (\textrm{there exists } X \HukuDiff \mathbb{E}X) = \mathbb{P} \left(\omega \ge \frac{1}{2}\right) = \frac{1}{2},
\]
and hence $\mathbb{E}X \not\in \huku{X}$.
\\
Actually we can show that $\huku{X}=\{\indicator{0}\}$ and hence $\hukuMin{X}=\indicator{0}$. In fact, by absurd let $B\in\huku{X}$ with $B\neq \indicator{0}$, then there exists $\alpha\in [0,1]$ such that $B_\alpha = [a,b]$ with $a<b$ and there exists $X_\alpha \HukuDiff B_\alpha$, here $\HukuDiff$ is considered as the Hukuhara difference for subsets in $\mathbb{R}$. On the other hand
\[
[-\omega, \omega] \HukuDiff [a , b] =
\left\{
\begin{array}{ll}
[-\omega -a, \omega -b], & \omega -b > -\omega-a,
\\
\{ -\frac{b+a}{2} \}, & \omega = \frac{b-a}{2},
\\
\textrm{it does not exist}, & \omega < \frac{b-a}{2},
\end{array}
\right.
\]
and, as consequence,
\[
\mathbb{P} (  [-\omega, \omega] \HukuDiff [a , b]  \textrm{ does not exist} ) = \misura \left[ \left(-\infty, \frac{b-a}{2}\right) \cap [0,1] \right] > 0
\]
where the last inequality is due to the fact that, by hypothesis, ${b-a}>0$. This is an absurd since $B\in\huku{X}$ by hypothesis. Thus $\huku{X}=\{ \indicator{0} \}\neq \mathbb{E} X = \indicator{ [-\frac{1}{2} , \frac{1}{2}] }$.
\\
Finally, in order to produce a more general example, let us consider
\[
X= \indicator{[-\omega, \omega]} \Mink \indicator{ [-\frac{1}{2} , \frac{1}{2}] } = \indicator{ [-\omega -\frac{1}{2} , \omega + \frac{1}{2}] }
\]
so that, from Corollary~\ref{cor:X'=X+D_then_C'=C+D}, we immediately obtain that
\[
\indicator{[-1,1]} = \mathbb{E} X \neq \hukuMin{X} = \indicator{[-\frac{1}{2},\frac{1}{2}]}.
\]
Note that, this is a case in which $\hukuMin{X}$ is different from $\indicator{0}$.
\end{oss}

\section{Conclusion}

In this paper, we have proven that any square integrable \frv{} can be decomposed as $X=\hukuMin{X}\Mink Y$, where $\hukuMin{X}$ is a unique deterministic fuzzy convex compact set (i.e. in $\Fuzzy[kc]{\Banach}$) and $Y$ is an element of $L^2[\Omega;\Fuzzy[kc]{\Banach}]$. This decomposition leads us to characterize \frv{} translations for which $\hukuMin{X}=\mathbb{E} X$, where the expectation is in the Aumann sense.
\\
This fact is important, for example, in view of Proposition~\ref{pro:degenere_implies_translation} that allows us to defuzzificate the randomness of a \frv{} process $\{X_t\}_{t\ge 0}$ for which $\mathbb{E}X_t=\indicator{c}$ holds at any time $t$. In fact, since $X_t$ is a translation at each $t$, it can be interpreted simply as a random element on $\Banach$.
\\
In general, working with a centered $X$ in $L^2[\Omega;\Fuzzy[kc]{\Banach}]$ one may distinguish different cases:
\begin{itemize}
\item
the case of defuzzificated randomness, for which $\mathbb{E}X=\hukuMin{X}$ and hence $X$ is a translation.

\item
the case for which $\mathbb{E}X\not\in \huku{X}$ and $\hukuMin{X}=\indicator{0}$; the randomness of $X$ is totally fuzzy.

\item
the case for which $\mathbb{E}X\not\in \huku{X}$ and $\hukuMin{X}\neq \indicator{0}$. In this case, one may take advantages from decomposition $X=\hukuMin{X} \Mink Y$ splitting the deterministic case from the random one.
\end{itemize}

Our decomposition of $\hukuMin{X}$ is a particular case where the problem posed in \cite[p.174--175]{mol05} is well--solved by defining the Hukuhara set $\huku{X}$. In this view, $\hukuMin{X}$ may be interpreted as an expectation for $X$ that satisfies some of the properties, listed in \cite[p.190]{mol05} for random closed sets but trivially extendible in the fuzzy case, of a \lq\lq reasonable\rq\rq{} expectation of $X$.
%
%

The decomposition theorem proposed in Section~\ref{sec:hukuhara_decomposition} could not be compared with the fuzzy regression problem stated in \cite{wun:nat02}. In fact, in that paper, the authors look for the best linear approximation function of a given square integrable \frv{} $Y$ by another square integrable \frv{} $X$, studying the minimization problem
\[
\inf_{ a\in\mathbb{R}, B\in\Fuzzy[kc]{\Banach} } \mathbb{E} [ d_2(Y, aX\Mink B)^2].
\]

Future works may consider the possibility to relax some hypothesis; for example, replacing $\mathbb{R}^d$ with an Hilbert or a Banach space (problems may arise considering the embedding $j$ and hence the closure of the Hukuhara set $\huku{X}$), or dropping convexity hypothesis and hence stating a decomposition theorem for a fuzzy random element whose level sets are not necessarily convex. Finally, note that we restricted our studies to the existence of a such $\hukuMin{X}$; however, it is certainly interesting to establish whenever $\hukuMin{X}$ could be explicitly computed, though even in particular cases.

\section*{Acknowledgements}

The authors would like to thank Prof.~V.~Capasso for the helpful, fruitful and long discussions during the preparation of this paper.

\end{document}